\documentclass{amsart}\usepackage{amsmath,amssymb}\usepackage[all]{xy}%

\newtheorem{theorem}{Theorem}
\newcommand{\bt}{\begin{theorem}}
\newcommand{\et}{\end{theorem}}

\newtheorem{lemma}{Lemma}
\newcommand{\bl}{\begin{lemma}}
\newcommand{\el}{\end{lemma}}

\newtheorem{corollary}{Corollary}
\newcommand{\bc}{\begin{corollary}}
\newcommand{\ec}{\end{corollary}}

\newcommand{\beq}{\begin{equation}}
\newcommand{\eeq}{\end{equation}}
\newcommand{\benum}{\begin{enumerate}}
\newcommand{\eenum}{\end{enumerate}}

\newcommand{\N}{\ensuremath{ \mathbf N }}

\newcommand{\mbF}{\ensuremath{\mathbf F}}

\newcommand{\mci}{\ensuremath{ \mathcal I}}

\newcommand{\mfI}{\ensuremath{ \mathfrak I}}

\DeclareMathOperator{\kernel}{\text{kernel}}

\date{\today}\author{Melvyn B. Nathanson}\address{Department of Mathematics\\Lehman College (CUNY)\\Bronx, NY 10468} \email{melvyn.nathanson@lehman.cuny.edu}\title{Dimensions of monomial varieties}
\begin{document}\maketitle

\begin{abstract}  The dimensions of certain varieties defined by monomials are computed using only high school algebra.  \end{abstract}

\section{Krull dimension and varieties}    \label{Krull:section:variety}

In this paper, a ring $R$ is a commutative ring with a multiplicative identity, and 
a field \mbF\  is an infinite field of any characteristic.

Let  $S$ be a nonempty set of polynomials 
in ${\mbF}[t_1,\ldots, t_n]$.
The \emph{ variety}\index{variety} (also called the \emph{algebraic set}) 
$V$ determined by $S$ 
is the set of points  in ${\mbF}^n$ 
that are common zeros of the polynomials in $S$, that is,
\[
V = V(S) = \left\{ (x_1 , \ldots, x_n ) \in {\mbF}^n: f(x_1,\ldots, x_n) = 0 \text{ for all } f \in S \right\}.
\]
The \emph{vanishing ideal}\index{vanishing ideal}  ${\mfI}(V)$ 
is the set of polynomials that vanish on the  variety $V$, that is, 
\[
{\mfI}(V) =   \left\{ f \in{\mbF}[t_1,\ldots, t_n] : f(x_1,\ldots, x_n) = 0 \text{ for all }  (x_1 , \ldots, x_n ) \in V \right\}.
\]
We have $S \subseteq \mfI(V)$, and so $\mfI(V)$ contains the ideal generated by $S$.
The quotient ring 
\[
{\mbF}(V) ={\mbF}[t_1,\ldots, t_n]/{\mfI}(V)
\]
is called the \emph{coordinate ring}\index{coordinate ring} of $V$.

Note that $\mfI(V)$ contains $S$, and so $\mfI(V)$ contains the ideal generated by $S$.

A \emph{prime ideal chain of length $n$}\index{chain!prime ideal}\index{prime ideal chain} 
in the ring $R$ is a strictly increasing sequence of $n+1$ prime ideals of $R$.  
The  \emph{Krull dimension}\index{Krull dimension} of $R$  
is the supremum of the lengths of  prime ideal chains in $R$.   
We define the dimension of the variety $V$ as the Krull dimension of 
its coordinate ring ${\mbF}(V)$.  

It is a basic theorem in commutative algebra  that  
the polynomial ring $\mbF[t_1,\ldots, t_n]$ 
has Krull dimension $n$.     
(Nathanson~\cite{nath17y} gives an elementary proof.  
Other references are  Atiyah and Macdonald~\cite[Chapter 11]{atiy-macd69}, 
Cox, Little, and O'Shea~\cite[Chapter 9]{cox-litt-oshe07}, and 
Kunz~\cite[Chapter 2]{kunz13}).  
If $S = \{0\} \subseteq \mbF[t_1,\ldots, t_n]$ is  the set whose only element 
is the zero polynomial, then $V = V(  \{0\} ) = \mbF^n$, and the vanishing ideal of $V$ 
is $\mci(V) = \mci(\mbF^n) =  \{0\}$.  
We obtain the coordinate ring 
\[
{\mbF}(V) = {\mbF}[t_1,\ldots, t_n]/{\mfI}(V) \cong {\mbF}[t_1,\ldots, t_n]
\]
and so the variety $\mbF^n$ has dimension $n$.

We adopt standard polynomial notation.  Let $\N_0$ denote the set of nonnegative integers.  
Associated to every $n$-tuple $I = (i_1,\ldots, i_n) \in \N_0^n$ is the monomial 
\[
t^I = t_1^{i_1} \cdots t_n^{i_n}.
\]
Every polynomial $f \in R[t_1,\dots, t_n]$ can be represented uniquely in the form 
\[
f = \sum_{I  \in \N_0^n} c_I t^I
\]
where $c_I \in R$ and $c_I \neq 0$ for only finitely many $I  \in \N_0^n$.

In this paper, two results about polynomials from high school algebra 
will enable us to compute the dimensions of certain  varieties defined by monomials.  
The first result is a factorization formula, and the second result follows 
from the fact that a polynomial of degree $d$ has at most $d$ roots in a field.

\bl                  \label{Krull:lemma:factorizationFormula}
For every nonnegative integer $i$, there is the polynomial identity 
\beq         \label{Krull:factorizationFormula}
u^i - v^i = (u-v) \Delta_i(u,v).
\eeq
where 
\beq         \label{Krull:factorization}
\Delta_i(u,v) = \sum_{j =0}^{i-1} u^{i-1-j} v^j.
\eeq
\el

\begin{proof}
We have 
\begin{align*}
(u-v) \sum_{j =0}^{i-1} u^{i-1-j} v^j 
& =  \sum_{j =0}^{i-1} u^{i-j} v^j - \sum_{j =0}^{i-1} u^{i-1-j} v^{j+1} \\
& =  \sum_{j =0}^{i-1} u^{i-j} v^j - \sum_{j =1}^{i} u^{i-j} v^{j} \\
& = u^i - v^i.  
\end{align*}
\end{proof}

\bl                      \label{Krull:lemma:polynomialZero}
Let \mbF\ be an infinite field.  
A polynomial $f \in \mbF[t_1,\dots, t_n]$ satisfies $f(x_1,\dots, x_n) = 0$ 
for all $(x_1,\dots, x_n) \in \mbF^n$ if and only if $f = 0$.  
\el

\begin{proof}
The proof is by induction on $n$.  
Let $n = 1$.  A nonzero polynomial $f \in \mbF[t_1]$ of degree $d$ 
has at most $d$ roots in \mbF, and so $f(x_1) \neq 0$ for some $x_1\in \mbF$.  
Thus, if  $f(x_1) = 0$ for all$x_1\in \mbF$, then $f = 0$. 

Let $n \geq 1$, and assume that the Lemma holds for polynomials in $n$ variables.
Let $f \in  \mbF[t_1,\dots, t_n,  t_{n+1}]$ have degree $d$ in the variable $ t_{n+1}$.  
There exist polynomials $f_i \in  \mbF[t_1,\dots, t_n]$ such that 
\[
f = f(t_1,\dots, t_n,  t_{n+1})  = \sum_{i=0}^d f_i(t_1,\dots, t_n) t_{n+1}^i
\]
For all $(x_1,\ldots, x_n) \in \mbF^n$, the polynomial  
\[
g(t_{n+1}) = f(x_1,\ldots, x_n, t_{n+1})  = \sum_{i=0}^d f_i(x_1,\dots, x_n) t_{n+1}^i 
\in \mbF[t_{n+1}]
\]
satisfies $g( x_{n+1}) = 0$ for all $x_{n+1} \in \mbF$, and so $g=0$.  Therefore, 
$ f_i(x_1,\dots, x_n) = 0$ for all $(x_1,\dots, x_n)  \in \mbF^n$ and $i=0,1,\ldots, d$.  
By the induction hypothesis, $f_i = 0$ for all $i=0,1,\ldots, d$, and so $f =0$.
\end{proof}

\section{An example of a plane curve}. \label{Krull:section:curve}

A \emph{hypersurface} is a variety that is the set of zeros of one nonzero polynomial.  
A \emph{plane algebraic curve} is a hypersurface in $\mbF^2$.
In this section we compute the dimension 
of a hypersurface $V$ in $\mbF^{m+1}$ defined by a polynomial of the form 
\[
f^* = t_{m+1} -  \lambda t_1^{a_1}t_2^{a_2}\cdots t_m^{a_m}  
\]
where $\lambda \in \mbF$ and $(a_1,\ldots, a_m) \in \N_0^m$.   
We shall prove that the \emph{monomial hypersurface}
\begin{align*}
V 
& = \{ (x_1,\ldots, x_m, x_{m+1}) \in \mbF^{m+1}: 
f^* (x_1,\ldots, x_m, x_{m+1}) = 0 \} \\ 
& = \{ (x_1,\ldots, x_m, x_{m+1}) \in \mbF^{m+1}: 
x_{m+1} = \lambda x_1^{a_1}x_2^{a_2}\cdots x_m^{a_m} \}
\end{align*}
has dimension $m$.

We begin with an example.
Consider the monomial $4t_1^3$ and the curve in $\mbF^2$ 
defined by the polynomial 
\[
f^* = t_2 - 4t_1^3 \in \mbF[t_1,t_2].
\]
Let
\[
V = \left\{ (x_1,x_2) \in \mbF^2: f^*(x_1,x_2) = 0  \right\} 
 = \left\{ (x_1,x_2) \in \mbF^2: x_2 = 4x_1^3  \right\}. 
\]
We shall prove that the vanishing ideal $\mfI(V)$ is the principal ideal
 generated by $f^*$.
 Because $f^* \in \mfI(V)$, it suffices to show that every 
 polynomial in $\mfI(V)$ is divisible by $f^*$.  
 
For $I = (i_1,i_2) \in \N_0^2$ and $t^I  = t_1^{i_1} t_2^{i_2}$, let 
\[
b_1 = i_1 + 3i_2
\]
and let $\Delta_{i_2}$ be the polynomial defined  by~\eqref{Krull:factorization} 
in Lemma~\ref{Krull:lemma:factorizationFormula}.  
We have
\begin{align*}
t^I - 4^{i_2} t_1^{b_1} 
& = t_1^{i_1} t_2^{i_2} - 4^{i_2} t_1^{i_1 + 3i_2} 
 = t_1^{i_1} \left( t_2^{i_2} - 4^{i_2} t_1^{ 3i_2} \right) \\
& = t_1^{i_1} \left( t_2^{i_2} - \left( 4 t_1^{ 3} \right)^{i_2} \right) 
 =  t_1^{i_1} \Delta_{i_2} (  t_2, 4 t_1^{ 3} ) \left( t_2 - 4 t_1^{ 3}  \right) \\
& = g_I f^*
\end{align*}
where $g_I =   t_1^{i_1} \Delta_{i_2} (  t_2, 4 t_1^{ 3} ) \in \mbF[t_1,t_2]$.  

Every polynomial $f  \in \mbF[t_1,t_2]$ 
can be represented uniquely in the form 
\[
f  = \sum_{ I = (i_1, t_2) \in \N_0^2} c_I t_1^{i_1} t_2^{i_2} 
= \sum_{ b_1 \in \N_0 } \quad
\sum_{ \substack{ I = (i_1, i_2) \in \N_0^2   
\\ i_1+  3 i_2= b_1 } }
c_I   t_1^{i_1}t^{i_2} 
\]
A polynomial $f \in \mbF[t_1,\ t_2]$ is in the vanishing ideal 
$\mfI(V)$ if and only if, for all  $x_1 \in \mbF$, 
\begin{align*}
0 & = f (x_1, 4 x_1^3 ) 
=   \sum_{ b_1\in \N_0 }    \quad
\sum_{ \substack{ I = (i_1, i_2) \in \N_0^2  
\\ i_1 + 3 i_2 = b_1}}  c_I 
 x_1^{i_1 } \left( 4 x_1^3 \right)^{i_2}  \\
 & =   \sum_{  b_1\in \N_0}    \quad
\sum_{ \substack{ I = (i_1, i_2) \in \N_0^2  \\ i_1 + 3 i_2 = b_1  }}
 c_I  4^{i_2} x_1^{i_1+3i_{m+1} }  \\
 & =  \sum_{  b_1\in \N_0}    \left(
\sum_{ \substack{ I = (i_1,i_2) \in \N_0^2 
\\ i_1 + 3i_2= b_1   }}
 c_I 4^{i_2} \right) x_1^{b_1}.  
\end{align*}  
By Lemma~\ref{Krull:lemma:polynomialZero}, 
because \mbF\ is an infinite field, the coefficients of this polynomial 
are zero, and so  
\[
\sum_{ \substack{ I = (i_1,  i_2) \in \N_0^2  
\\ i_1 + 3 i_2= b_1 }}
 c_I 4^{i_2}  = 0
\]
for all $b_1 \in \N_0$.  
The ordered pair $I =  (b_1, 0)$ is one of the terms in this sum, and so 
\[
-c_{(b_1,0)} =  \sum_{ \substack{ I = (i_1,  i_) \in \N_0^2,   
\\ i_{\ell} + 3 i_{m+1} = b_1  \\  I \neq  (b_1,0)}}
 c_I 4^{i_2}.  
\]
Therefore, $f \in \mfI(V)$ implies 
\begin{align*}
f 
& =  \sum_{  b_1\in \N_0}    \left( c_{(b_1,0)} t_1^{b_1} + 
\sum_{ \substack{ I = (i_1,i_2) \in \N_0^2 
\\ i_1 + 3i_2= b_1 \\ I \neq (b_1,0)  }}
 c_I t_1^{i_1} t_2^{i_2} \right)  \\
 & =  \sum_{  b_1\in \N_0}    \left( 
\sum_{ \substack{ I = (i_1,i_2) \in \N_0^2 
\\ i_1 + 3i_2= b_1 \\ I \neq (b_1,0)  }}
 c_I t_1^{i_1} t_2^{i_2} -  \sum_{ \substack{ I = (i_1,  i_) \in \N_0^2,   
\\ i_{\ell} + 3 i_{m+1} = b_1  \\  I \neq  (b_1,0)}}
 c_I 4^{i_2}   t_1^{b_1} \right)  \\
  & =  \sum_{  b_1\in \N_0}   \quad 
\sum_{ \substack{ I = (i_1,i_2) \in \N_0^2 
\\ i_1 + 3i_2= b_1 \\ I \neq (b_1,0)  }}
 c_I   \left(  t_1^{i_1} t_2^{i_2} - 4^{i_2}   t_1^{b_1} \right)  \\
 & =  \sum_{  b_1\in \N_0}   \quad 
\sum_{ \substack{ I = (i_1,i_2) \in \N_0^2 
\\ i_1 + 3i_2= b_1 \\ I \neq (b_1,0)  }}
 c_I  g^I f^* 
 \end{align*}
and so $f^*$ divides $f$.  Thus, every polynomial $f \in \mfI(V)$ 
is contained in the principal ideal generated by $f^*$.

The function 
\[
\varphi:\mbF[t_1, t_2] \rightarrow \mbF[t_1] 
\]
defined by 
\[
\varphi(t_1) = t_1
\]
and
\[
\varphi(t_2) = 4 t_1^3
\]
is a surjective ring homomorphism with 
\[
\kernel(\varphi) = \left\{ f \in  \mbF[t_1, t_2] : 
f(t_1, 4^{i_2}t_1^3 )= 0 \right\} 
= \mfI(V).
\]
Therefore, 
\[
\mbF[V]= \mbF[t_1, t_2]/\mci (V)\cong \mbF[t_1].
\]
The polynomial ring $\mbF[t_1]$ has Krull dimension $1$, 
and so the coordinate ring $\mbF(V)$ of the curve has Krull dimension $1$ 
and the curve  has  dimension $1$.

\section{Dimension of a monomial hypersurface}
We shall prove that every monomial hypersurface in $\mbF^{m+1}$ 
has dimension $m$.  
The proof is elementary, like the proof in Section~\ref{Krull:section:curve}, 
but a bit more technical.

\bl                   \label{Krull:lemma:reduce-1}
For $\lambda \in \mbF$ and $(a_1,\ldots, a_m) \in \N_0^m$, 
consider the polynomial 
\[
f^* =  t_{m+1} -  \lambda t_1^{a_1} \cdots t_m^{a_m} 
\in \mbF[t_1,\ldots, t_{m+1}].
\]
For $I = (i_1,\ldots, i_m,i_{m+1}) \in \N_0^{m+1}$, let 
\[
 b_{\ell}  = i_{\ell} +  a_{\ell} i_{m+1} \qquad \text{for $\ell = 1, \ldots, m$}.  
\]
There exists a polynomial $g_I \in \mbF[t_1,\ldots, t_{m+1}]$ such that 
\[
 t^I  - \lambda^{i_{m+1}}t_1^{b_1} \cdots t_m^{b_m} = g_I f^*.
\]
\el

\begin{proof}
Let  $\Delta_i(u,v)$ be the polynomial defined by~\eqref{Krull:factorization}.
We have 
\begin{align*}
t^I  - & \lambda^{i_{m+1}}t_1^{b_1} \cdots t_m^{b_m} \\
& = t_1^{i_1} \cdots t_m^{i_m} t_{m+1}^{i_{m+1}} 
-  \lambda^{i_{m+1}}t_1^{i_1+a_1i_{m+1}} \cdots t_m^{i_m+a_m i_{m+1}} \\
& =  t_1^{i_1} \cdots t_m^{i_m} \left(  t_{m+1}^{i_{m+1}} 
-  \lambda^{i_{m+1}}t_1^{a_1i_{m+1}} \cdots t_m^{a_m i_{m+1}} \right) \\
& =  t_1^{i_1} \cdots t_m^{i_m} \left(  t_{m+1}^{i_{m+1}} 
- \left( \lambda t_1^{a_1} \cdots t_m^{a_m} \right)^{i_{m+1}} \right) \\
& = t_1^{i_1} \cdots t_m^{i_m} \Delta_{i_{m+1}}(t_{m+1}, \lambda t_1^{a_1} \cdots t_m^{a_m} ) 
\left( t_{m+1} -  \lambda t_1^{a_1} \cdots t_m^{a_m} \right) \\
& = g_I f^*
\end{align*}
where 
\[
g_I =  t_1^{i_1} \cdots t_m^{i_m} \Delta_{i_{m+1}}(t_{m+1}, \lambda t_1^{a_1} \cdots t_m^{a_m} ) 
\in \mbF[t_1,\ldots, t_{m+1}].  
\]  
This completes the proof.  
\end{proof}

\bt
Let \mbF\ be an infinite field.  
For $\lambda \in \mbF$ and $(a_1,\ldots, a_m) \in \N_0^m$, 
consider the polynomial 
\[
f^* =  t_{m+1}  - \lambda t_1^{a_1}t_2^{a_2}\cdots t_m^{a_m}  
\in \mbF[t_1,\ldots, t_m, t_{m+1}]
\]  
and the associated hypersurface 
\begin{align*}
V 
& =  \left\{ \left(x_1,\ldots, x_m, x_{m+1} \right) 
\in \mbF^{m+1} : f^*\left(x_1,\ldots, x_m, x_{m+1} \right) = 0  \right\}. \\
& = \left\{ \left(x_1,\ldots, x_m,   \lambda x_1^{a_1}x_2^{a_2}\cdots x_m^{a_m} \right) 
\in \mbF^{m+1} :  (x_1,\ldots, x_m)   \in \mbF^m \right\}.
\end{align*}
The vanishing ideal $\mfI(V)$ is the principal ideal generated by $f^*$.  
\et

\begin{proof}
The vanishing ideal $\mfI(V)$ contains $f^*$, and so 
$\mfI(V)$ contains the principal ideal generated by $f^*$. 
Therefore, it suffices to prove that $\mfI(V)$ 
is contained in the principal ideal generated by $f^*$. 

For every $(m+1)$-tuple $ I = (i_1, \ldots, i_m, i_{m+1}) \in \N_0^{m+1}$, 
there is a unique $m$-tuple $(b_1,\ldots, b_m) \in \N_0^m$ such that 
\[
i_{\ell} + a_{\ell} i_{m+1} = b_{\ell}
\]
for $\ell = 1, \ldots, m$.  
Thus, every polynomial $f  \in \mbF[t_1,\ldots, t_m, t_{m+1}]$ 
can be represented uniquely in the form 
\begin{align*}
f  & = \sum_{I  \in \N_0^{m+1} } c_I t^I \\
& = \sum_{ (b_1,\ldots, b_m) \in \N_0^m }    \quad
\sum_{ \substack{ I = (i_1, \ldots, i_m, i_{m+1}) \in \N_0^{m+1}   
\\ i_{\ell} + a_{ \ell} i_{m+1} = b_{\ell}  \\  \text{for $\ell = 1,\ldots, m$}   
   }}
c_I   t_1^{i_1} \cdots t_m^{i_m} t_{m+1}^{i_{m+1}}.
\end{align*}
A polynomial $f \in \mbF[t_1,\ldots, t_m, t_{m+1}]$ is in the vanishing ideal 
$\mfI(V)$ if and only if, for all  $(x_1,\ldots, x_m) \in \mbF^m $, 
\begin{align*}
0 & = f (x_1,\ldots, x_m,   \lambda x_1^{a_{1}} x_2^{a_{2}}\cdots x_m^{a_{m}}  ) \\
& =   \sum_{ (b_1,\ldots, b_m) \in \N_0^m }    \quad
\sum_{ \substack{ I = (i_1, \ldots, i_m, i_{m+1}) \in \N_0^{m+1}   
\\ i_{\ell} + a_{ \ell} i_{m+1} = b_{\ell}  \\  \text{for $\ell = 1,\ldots,m$}   
   }}
 c_I 
 x_1^{i_1 } \cdots  x_m^{i_m }  
 \left( \lambda x_1^{a_{1}} x_2^{a_{2}}\cdots x_m^{a_{m}} \right)^{i_{m+1}}  \\
 & =   \sum_{ (b_1,\ldots, b_m) \in \N_0^m }    \quad
\sum_{ \substack{ I = (i_1, \ldots, i_m, i_{m+1}) \in \N_0^{m+1}   
\\ i_{\ell} + a_{ \ell} i_{m+1} = b_{\ell}  \\  \text{for $\ell = 1,\ldots,m$}   
   }}
 c_I   \lambda^{i_{m+1}} x_1^{i_1+a_1i_{m+1} } \cdots  x_m^{i_m +a_m i_{m+1} }  \\
 & =  \sum_{ (b_1,\ldots, b_m) \in \N_0^m }    \left(
\sum_{ \substack{ I = (i_1, \ldots, i_m, i_{m+1}) \in \N_0^{m+1}   
\\ i_{\ell} + a_{ \ell} i_{m+1} = b_{\ell}  \\  \text{for $\ell = 1,\ldots,m$}   
   }}
 c_I \lambda^{i_{m+1}} \right) x_1^{b_1}\cdots x_m^{b_m}.
\end{align*}  
By Lemma~\ref{Krull:lemma:polynomialZero}, 
the coefficients of this polynomial are zero, and so  
\beq                       \label{Krull:b-sum}
\sum_{ \substack{ I = (i_1, \ldots, i_m, i_{m+1}) \in \N_0^{m+1}   
\\ i_{\ell} + a_{ \ell} i_{m+1} = b_{\ell}  \\  \text{for $\ell = 1,\ldots,m$}   
   }}
 c_I \lambda^{i_{m+1}}  = 0
\eeq
for all $ (b_1,\ldots, b_m) \in \N_0^m $.  
The $(m+1)$-tuple $I =  (b_1,\ldots, b_m,0)$ is one of the terms 
in the sum~\eqref{Krull:b-sum}, and so 
\[
-c_{(b_1,\ldots, b_m,0)} =  \sum_{ \substack{ I = (i_1, \ldots, i_m, i_{m+1}) \in \N_0^{m+1},   
\\ i_{\ell} + a_{ \ell} i_{m+1} = b_{\ell}  \\  \text{for $\ell = 1,\ldots,m$,}   
\\  I \neq  (b_1,\ldots, b_m,0)}}
 c_I \lambda^{i_{m+1}}.  
\]
Therefore, $f \in \mfI(V)$ implies 
\begin{align*}
f 
&  = \sum_{ (b_1,\ldots, b_m) \in \N_0^m }    \quad
\sum_{ \substack{ I = (i_1, \ldots, i_m, i_{m+1}) \in \N_0^{m+1}   
\\ i_{\ell} + a_{ \ell} i_{m+1} = b_{\ell}  \\  \text{for $\ell = 1,\ldots,m$}  }}
c_I    t^I \\
 &  = \sum_{ (b_1,\ldots, b_m) \in \N_0^m }    \left(
\sum_{ \substack{ I = (i_1, \ldots, i_m, i_{m+1}) \in \N_0^{m+1}   
\\ i_{\ell} + a_{ \ell} i_{m+1} = b_{\ell}  \\  \text{for $\ell = 1,\ldots,m$} \\ 
 I \neq  (b_1,\ldots, b_m,0) }}
c_I   t^I  + c_{(b_1,\ldots, b_m,0)} t_1^{b_1} \cdots t_m^{b_m} 
\right)\\
 &  = \sum_{ (b_1,\ldots, b_m) \in \N_0^m }    \left(
\sum_{ \substack{ I = (i_1, \ldots, i_m, i_{m+1}) \in \N_0^{m+1}   
\\ i_{\ell} + a_{ \ell} i_{m+1} = b_{\ell}  \\  \text{for $\ell = 1,\ldots,m$} \\ 
 I \neq  (b_1,\ldots, b_m,0)}}
c_I   t^I  - \sum_{ \substack{ I = (i_1, \ldots, i_m, i_{m+1}) \in \N_0^{m+1}   
\\ i_{\ell} + a_{ \ell} i_{m+1} = b_{\ell}  \\  \text{for $\ell = 1,\ldots,m$}   
\\ I \neq  (b_1,\ldots, b_m,0)}}
 c_I \lambda^{i_{m+1}} t_1^{b_1} \cdots t_m^{b_m} 
\right)\\
 &  = \sum_{ (b_1,\ldots, b_m) \in \N_0^m }   \quad 
\sum_{ \substack{ I = (i_1, \ldots, i_m, i_{m+1}) \in \N_0^{m+1}   
\\ i_{\ell} + a_{ \ell} i_{m+1} = b_{\ell}  \\  \text{for $\ell = 1,\ldots, m$} \\ 
 I \neq  (b_1,\ldots, b_m,0)} }
c_I  \left(  t^I  - \lambda^{i_{m+1}}t_1^{b_1} \cdots t_m^{b_m} 
\right)\\
& = \sum_{ (b_1,\ldots, b_m) \in \N_0^m }   \quad 
\sum_{ \substack{ I = (i_1, \ldots, i_m, i_{m+1}) \in \N_0^{m+1}   
\\ i_{\ell} + a_{ \ell} i_{m+1} = b_{\ell}  \\  \text{for $\ell = 1,\ldots, m$} \\  I \neq  (b_1,\ldots, b_m,0)} }
c_I g_I f^*
\end{align*}
by Lemma~\ref{Krull:lemma:reduce-1},  
and so $f$ is in the principal ideal generated by $f^*$.
This completes the proof.  
\end{proof}

\bt
For $\lambda \in \mbF$ and $(a_1,\ldots, a_m) \in \N_0^m$, the hypersurface 
\[
V =  \left\{ \left(x_1,\ldots, x_m,   \lambda x_1^{a_1}x_2^{a_2}\cdots x_m^{a_m} \right) 
\in \mbF^{m+1} :  (x_1,\ldots, x_m)   \in \mbF^m \right\} 
\]
has dimension m.
\et

\begin{proof}
The function 
\[
\varphi:\mbF[t_1,\ldots, t_{m+1}] \rightarrow \mbF[t_1,\ldots, t_m] 
\]
defined by 
\[
\varphi(t_{\ell}) = t_{\ell}  \qquad \text{for $\ell = 1,\ldots, m$}
\]
and
\[
\varphi(t_{m+1}) = \lambda t_1^{a_1} \cdots t_m^{a_m} 
\]
is a surjective ring homomorphism with 
\[
\kernel(\varphi) = \left\{ f \in  \mbF[t_1,\ldots, t_m] : 
f(t_1,\ldots, t_m, \lambda^{i_{m+1}}t_1^{a_1} \cdots t_m^{a_m} = 0 \right\} 
= \mfI(V).
\]
Therefore, 
\[
\mbF[V]= \mbF[t_1,\ldots, t_m, t_{m+1}]/\mci (V)\cong \mbF[t_1,\ldots, t_m].
\]
The polynomial ring  $\mbF[t_1,\ldots, t_m]$ has Krull dimension $m$, and so 
the coordinate ring $\mbF[V]$ has Krull dimension $m$ 
and the hypersurface $V$ has dimension $m$.
This completes the proof.  
\end{proof}

\section{Varieties defined by several monomials}

Let $m$ and $k$ be positive integers, and let $n = m+k$.
For $j = 1,2,\ldots, k$, let $\lambda_j \in \mbF$ 
and $(a_{1,j}, a_{2,j} ,\ldots, a_{m,j}) \in \N_0^m$.
Consider the polynomials
\beq           \label{Krull:define_f}
f_j^* =  t_{m+j}  - \lambda_j t_1^{a_{1,j}}t_2^{a_{2,j}}\cdots t_m^{a_{m,j}} 
 \in  \mbF[t_1,\ldots,t_n].
\eeq
Let $V$ be the  variety in $\mbF^n$ determined by the set of polynomials 
\[
S = \{f^*_j:j=1,\ldots, k\} 
\]
and let $\mfI(V)$ be the vanishing ideal of $V$.  
We shall prove that the coordinate ring 
$\mbF[V] = \mbF[t_1,\ldots,t_n]/\mfI(V)$ is isomorphic 
to the polynomial ring $\mbF[t_1,\ldots,t_m]$, and so $V$ has dimension $m$.

\bl                   \label{Krull:lemma:reduce-2}
Let $R$ be a ring.  
For $j = 1,2,\ldots, k$, let $\lambda_j \in \mbF$ 
and $(a_{1,j}, a_{2,j} ,\ldots, a_{m,j}) \in \N_0^m$. 
Define the polynomial  $f_j^* $ by~\eqref{Krull:define_f}.  
For $I = (i_1,\ldots, i_m, i_{m+1}, \ldots, i_{m+k}) \in \N_0^{m+k}$, let 
\[
 b_{\ell}  = i_{\ell} +  \sum_{j=1}^k a_{\ell, j} i_{m+j} \qquad \text{for $\ell = 1, \ldots, m$}.  
\]  
There exist polynomials $g_{I,1}, \ldots, g_{I,k} \in \mbF[t_1,\ldots, t_{m+k}]$ such that 
\beq                          \label{Krull:reduce-2}
t^I - \prod_{j=1}^k \lambda_j^{i_{m+j}}  \ t_1^{b_1} \cdots t_m^{b_m} = \sum_{j=1}^k g_{I,j} \ f_j^*.
\eeq
\el

\begin{proof}
The proof is by induction on $k$.  The case $k=1$ is Lemma~\ref{Krull:lemma:reduce-1}.  
Assume that Lemma~\ref{Krull:lemma:reduce-2} is true for the positive integer $k$.  
We shall prove the Lemma  for $k+1$.

For
\[
I = (i_1,\ldots, i_{m+k}) \in \N_0^{m+k}
\]
and
\[
I'  = (i_1,\ldots, i_{m+k+1}) \in \N_0^{m+k+1} 
\]
we have 
\[
t^{I}  = \prod_{\ell=1}^{m} t_{\ell}^{i_{\ell}} \prod_{j=1}^{k}t_{m+j}^{i_{m+j}} 
\]
and 
\[
t^{I'}  = \prod_{\ell=1}^{m}t_{\ell}^{i_{\ell}}  \prod_{j=1}^{k+1}t_{m+j}^{i_{m+j}} 
= t^I \ t_{m+k+1}^{i_{m+k+1}}.
\]
For $\ell = 1,\ldots, m$, define 
\[
b_{\ell}  =  i_{\ell} +  \sum_{j=1}^{k} a_{\ell, j} i_{m+j}
\]
and
\[
b'_{\ell}  =  i_{\ell} +  \sum_{j=1}^{k+1} a_{\ell, j} i_{m+j} 
= b_{\ell} + a_{\ell, k+1} i_{m+ k+1}.
\]

We have 
\begin{align*}
t^{I'} &  -  \prod_{j=1}^{k+1} \lambda_j^{i_{m+j}}  \ t_1^{b'_1} \cdots t_m^{b'_m} \\
&  =  t_{m+k+1}^{i_{m+k+1}}  \left( t^I - \prod_{j=1}^k \lambda_j^{i_{m+j}}  \ t_1^{b_1} \cdots t_m^{b_m} \right) \\
& \qquad +  t_{m+k+1}^{i_{m+k+1}} \left( \prod_{j=1}^k \lambda_j^{i_{m+j}}  \ t_1^{b_1} \cdots t_m^{b_m} \right) -  \prod_{j=1}^{k+1} \lambda_j^{i_{m+j}}  \ t_1^{b'_1} \cdots t_m^{b'_m} 
\end{align*} 
By the induction hypothesis, 
there exist polynomials $g_{I,1}, \ldots, g_{I,k} \in \mbF[t_1,\ldots, t_{m+k}]$ 
that satisfy~\eqref{Krull:reduce-2}, and so 
\[
  t_{m+k+1}^{i_{m+k+1}} 
   \left( t^I - \prod_{j=1}^k \lambda_j^{i_{m+j}}  \ t_1^{b_1} \cdots t_m^{b_m} \right) 
=   t_{m+k+1}^{i_{m+k+1}}  \sum_{j=1}^k g_{I,j} \ f_j^*.
\]
Applying the factorization formula~\eqref{Krull:factorizationFormula}, we obtain 
\begin{align*}
 t_{m+k+1}^{i_{m+k+1}} & \left( \prod_{j=1}^k \lambda_j^{i_{m+j}}  \ t_1^{b_1} \cdots t_m^{b_m} \right) -  \prod_{j=1}^{k+1} \lambda_j^{i_{m+j}}  \ t_1^{b'_1} \cdots t_m^{b'_m} \\
& =  \left( \prod_{j=1}^k \lambda_j^{i_{m+j}}   \ t_1^{b_1} \cdots t_m^{b_m} \right)  
\left(  t_{m+k+1}^{i_{m+k+1}} - \lambda_{k+1}^{i_{m+k+1}} \prod_{\ell=1}^m  t_{\ell}^{ a_{\ell, k+1} i_{m+ k+1}} \right) \\
& =  \left( \prod_{j=1}^k \lambda_j^{i_{m+j}}   \ t_1^{b_1} \cdots t_m^{b_m} \right)  
\left(  t_{m+k+1}^{i_{m+k+1}} - \left( \lambda_{k+1} \prod_{\ell=1}^m  t_{\ell}^{ a_{\ell, k+1}}
\right)^{ i_{m+ k+1}}   \right) \\
& =  \left( \prod_{j=1}^k \lambda_j^{i_{m+j}}   \ t_1^{b_1} \cdots t_m^{b_m} \right)  
\Delta_{i_{m+k+1}}
 \left(  t_{m+k+1} , \  \lambda_{k+1} \prod_{\ell=1}^m  t_{\ell}^{ a_{\ell, k+1}} \right) 
  \left(  t_{m+k+1} -  \lambda_{k+1} \prod_{\ell=1}^m  t_{\ell}^{ a_{\ell, k+1}} \right) \\
& = g_{I,k+1} \ f^*_{k+1} 
\end{align*} 
where 
\[
 g_{I,k+1} = \left( \prod_{j=1}^k \lambda_j^{i_{m+j}}   \ t_1^{b_1} \cdots t_m^{b_m} \right)  
\Delta_{i_{m+k+1}}
 \left(  t_{m+k+1} , \     \lambda_{k+1} \prod_{\ell=1}^m  t_{\ell}^{ a_{\ell, k+1}} \right). 
\]
This completes the proof.  
\end{proof}

\bt                                                                 \label{Krull:theorem:poly-m}
Let $n = m+k$.  Let \mbF\ be an infinite field.   
For $j = 1,2,\ldots, k$, let $\lambda_j \in \mbF$ 
and $(a_{1,j}, a_{2,j} ,\ldots, a_{m,j}) \in \N_0^m$,  
and let $f_j^* $ be the polynomial defined by~\eqref{Krull:define_f}.  
Let $V \subseteq \mbF^n$ be the variety determined by the set 
$S = \{f_1^*,\ldots, f_k^*\} \subseteq \mbF[t_1,\ldots, t_{n}]$.  
The vanishing ideal $\mfI(V)$ is the ideal generated by $S$.  
\et

\begin{proof}
The ideal $\mfI(V)$ contains $S$, and so $\mfI(V)$ contains 
the ideal generated by $S$.  Thus, it suffices to prove that 
 ideal generated by $S$ contains every polynomial in $\mfI(V)$. 
 
The variety determined by $S$ is 
\[
V = \left\{ 
\left( x_1,\ldots, x_m, \lambda_1 x_1^{a_{1,1}} \cdots x_m^{a_{m,1}}, \ldots, 
\lambda_k x_1^{a_{1,k}} \cdots x_m^{a_{m,k}}  \right) :
(x_1,\ldots, x_m) \in \mbF^m \right\}.
\]
For every $(m+k)$-tuple $ I = (i_1, \ldots, i_m, i_{m+1},\ldots, i_{m+k}) \in \N_0^{m+k}$, 
there is a unique $m$-tuple $(b_1,\ldots, b_m) \in \N_0^m$ such that 
\[
i_{\ell} + \sum_{j=1}^k a_{\ell,j} i_{m+j} = b_{\ell}
\]
for  $\ell = 1, \ldots, m$.  Let
\[
\sum_{I  (b_1,\ldots, b_m)} = \quad \sum_{  \substack{ I = (i_1, \ldots, i_m, i_{m+1}, \ldots, i_{m+k} ) \in \N_0^{n}    \\  i_{\ell} +  \sum_{j=1}^k a_{\ell, j} i_{m+j} = b_{\ell}  \\  \text{for $\ell = 1,\ldots,m$}}}.
\]
Every polynomial $f \in \mbF[t_1,\ldots, t_n]$ has a unique representation in the form
\begin{align*}
f &  = \sum_{I \in \N_0^n} c_I t^I   
 = \sum_{ (b_1,\ldots, b_m) \in \N_0^m }  
 \sum_{I  (b_1,\ldots, b_m)} c_I   t^I 
\end{align*}
where $c_I \in \mbF$ and $c_I \neq 0$ for only finitely many $n$-tuples $I$.
If $f \in \mfI(V)$, then
\begin{align*}
0 
& = f \left( x_1,\ldots, x_m, \lambda_1 x_1^{a_{1,1}} \cdots x_m^{a_{m,1}}, \ldots, 
\lambda_k x_1^{a_{1,k}} \cdots x_m^{a_{m,k}}  \right) \\
& = \sum_{ (b_1,\ldots, b_m) \in \N_0^m }   \sum_{ I  (b_1,\ldots, b_m) } 
c_I   x_1^{i_1} \cdots  x_m^{i_m}
\left( \lambda_1 x_1^{a_{1,1}} \cdots x_m^{a_{m,1}} \right)^{i_{m+1}}
\cdots \left( \lambda_k x_1^{a_{1,k}} \cdots x_m^{a_{m,k}} \right)^{i_{m+k}} \\
& = \sum_{ (b_1,\ldots, b_m) \in \N_0^m }  
 \sum_{I  (b_1,\ldots, b_m)} c_I \prod_{j=1}^k \lambda_{j}^{i_{m+j}} 
\prod_{ \ell =1}^m x_{\ell}^{i_{\ell} + \sum_{j=1}^k a_{\ell,j} i_{m+j} }   \\
& = \sum_{ (b_1,\ldots, b_m) \in \N_0^m }   \left( 
 \sum_{I  (b_1,\ldots, b_m)} c_I \prod_{j=1}^k \lambda_{j}^{i_{m+j}} \right) 
x_1^{b_1} \cdots x_m^{b_m}
\end{align*}
for all $(x_1,\ldots, x_m) \in \mbF^m$.  
It follows from Lemma~\ref{Krull:lemma:polynomialZero} that  
\begin{align*}
0 & = \sum_{ I_{(b_1,\ldots, b_m)} } c_I \prod_{j=1}^k \lambda_{j}^{i_{m+j}} 
 = c_{(b_1,\ldots, b_m, 0, \ldots, 0) } +  
\sum_{ \substack{ I_{(b_1,\ldots, b_m)} \\ I \neq  (b_1,\ldots, b_m, 0, \ldots, 0) }}
c_I \prod_{j=1}^k \lambda_{j}^{i_{m+j}} 
\end{align*}
 for all $(b_1,\ldots, b_m) \in \N_0^m$, and so 
\begin{align*}
f 
& = \sum_{ (b_1,\ldots, b_m) \in \N_0^m }   \sum_{I  (b_1,\ldots, b_m)} c_I   t^I \\ 
& = \sum_{ (b_1,\ldots, b_m) \in \N_0^m }  
\left(  \sum_{ \substack{ I_{(b_1,\ldots, b_m)} \\ I \neq  (b_1,\ldots, b_m, 0, \ldots, 0) }} c_I   t^I
+ c_{ (b_1,\ldots, b_m, 0, \ldots, 0)} t_1^{b_1} \cdots t_m^{b_m} \right) \\
 & = \sum_{ (b_1,\ldots, b_m) \in \N_0^m }  
\left(  \sum_{ \substack{ I_{(b_1,\ldots, b_m)} \\ I \neq  (b_1,\ldots, b_m, 0, \ldots, 0) }} c_I   t^I
- \sum_{ \substack{ I_{(b_1,\ldots, b_m)} \\ I \neq  (b_1,\ldots, b_m, 0, \ldots, 0) }}
c_I \prod_{j=1}^k \lambda_{j}^{i_{m+j}} 
 t_1^{b_1} \cdots t_m^{b_m} \right) \\
 & = \sum_{ (b_1,\ldots, b_m) \in \N_0^m }  
\sum_{ \substack{ I_{(b_1,\ldots, b_m)} \\ I \neq  (b_1,\ldots, b_m, 0, \ldots, 0) }} 
c_I  \left(  t^I -  \prod_{j=1}^k \lambda_{j}^{i_{m+j}} 
 t_1^{b_1} \cdots t_m^{b_m} \right).
\end{align*}
Lemma~\ref{Krull:lemma:reduce-2} immediately implies that $f$ is in the ideal 
generated by $S$.
This completes the proof.  
\end{proof}

\bt                                                                 \label{Krull:theorem:dim-m}
The variety $V$ has dimension $m$.
\et

\begin{proof}
The function 
\[
\varphi:\mbF[t_1,\ldots, t_{m+1}, \ldots, t_{m+k}] \rightarrow \mbF[t_1,\ldots, t_m] 
\]
defined by 
\[
\varphi(t_{\ell}) = 
t_{\ell}  \qquad \text{for $\ell = 1,\ldots, m$}
\]
and
\[
\varphi(t_{m+j}) = \lambda_j   t_1^{a_{1,j}} \cdots t_m^{a_{m,j}} 
\qquad \text{for $j = 1,\ldots, k$}
\]
is a surjective ring homomorphism with 
\begin{align*}
&\kernel  (\varphi) \\
& = \left\{ f \in  \mbF[t_1,\ldots, t_n] : 
f \left(t_1,\ldots, t_m, \lambda_1t_1^{a_{1,1}} \cdots t_m^{a_{m,1}} ,
\ldots,  \lambda_k t_1^{a_{1,k}} \cdots t_m^{a_{m,k}}   \right) = 0 \right\} \\
& = \mfI(V).
\end{align*}
Therefore, 
\[
\mbF[V] = \mbF[t_1,\ldots, t_n]/\mci (V)\cong \mbF[t_1,\ldots, t_m]
\]
and the coordinate ring of $\mfI(V)$ has Krull dimension $m$.  
This completes the proof.  
\end{proof}

\def\cprime{$'$} \def\cprime{$'$} \def\cprime{$'$}
\providecommand{\bysame}{\leavevmode\hbox to3em{\hrulefill}\thinspace}
\providecommand{\MR}{\relax\ifhmode\unskip\space\fi MR }
\providecommand{\MRhref}[2]{%
  \href{http://www.ams.org/mathscinet-getitem?mr=#1}{#2}
}
\providecommand{\href}[2]{#2}

\end{document}